\documentclass[12pt]{article}

\begin{document}

\newtheorem{thm}{Theorem}
\newtheorem{prop}{Proposition}
\newtheorem{coroll}{Corollary}
\newtheorem{defn}{Definition}
\newtheorem{lemma}{Lemma}
\title{Infinitesimal aspects of the Laplace operator}
\author{Anders Kock}
\date{}

\maketitle

%{\bf Abstract.} 
In the context of synthetic differential 
geometry, we study the Laplace ope\-rator an a Riemannian 
manifold. The main new aspect is a neighbourhood of the 
diagonal, smaller than the  second neighbourhood usually 
required as support for second order differential operators. 
The new neighbourhood has the property that a function is affine on it 
if and only if it is harmonic.

 \section*{Introduction} Recall \cite{FM}, \cite{SDG}, \cite{Levi-Civita}, 
\cite{DFIC}  that any manifold $M$, when seen in a 
model of Synthetic Differential Geometry (SDG), carries a reflexive 
symmetric relation $\sim _k$ ($k=0,1,2,...$), where $x\sim _k y$ reads 
``$x$ and $y$ are $k$-neighbours''; $x \sim _0 y$ means $x=y$; $x\sim _k y 
$ implies $x\sim _{k+1} y$.  Also, $x\sim _k y , y\sim _l z$ implies 
$x\sim_{k+l} z$.  The set of $(x,y)\in 
M\times M$ with $x\sim _k y$ is denoted $M_{(k)}$, the ``$k$'th 
neighbourhood of the diagonal", and for fixed $x$, the set $\{y \in M 
\mid y\sim _k x \}$ is denoted ${\cal M}_k (x)$ (``the $k$-monad 
around $x$").  
In $R^n$, ${\cal M}_k (0)$ is denoted $D_k (n)$.  
Its elements $u$ are characterized by the condition that any 
homogeneous polynomial of degree $k+1$ vansihes on $u$.
--- In this context, a {\em Riemannian metric} on $M$ can 
be given in terms of a map
$$g: M_{(2)} \to R$$
with $g(x,x)=0$, and with $g(x,y)= g(y,x)$ to be thought of as 
the ``square-distance between $x$ and $y$", se 
\cite{Levi-Civita}, \cite{VF}. (Also $g$ should be 
positive-definite, in a certain sense.)

Given a Riemannian metric $g$ on $M$, in this sense, one can 
construct the Levi-Civita connection \cite{Levi-Civita}, 
volume form \cite{VF}, and hence also a notion of {\em 
divergence} of a vector field. And to a function $f:M \to R$, 
one can construct its {\em gradient vector field}, and hence 
one can construct the Laplacian $\Delta$ by $\Delta (f) = div 
(grad (f))$. This is what we shall {\em not} do here, rather, 
we shall exploit the richness of synthetic language to give a 
more economic and more geometric construction of $\Delta$. The 
construction is more economic in the sense that the definition 
of $\Delta f (x)$ only depends on knowing $f$ on a certain 
subset ${\cal M}_{L} (x) \subseteq {\cal M}_{2}(x)$, where ${\cal 
M}_{2}(x)$ is what is required to make the usual $div \; grad $ 
construction work, or for defining the individual terms in the formula 
$\Delta f (x) = \sum \partial ^2 f /\partial x_i ^2 (x)$.

The description of ${\cal M}_L  (x) \subseteq M$, or equivalently, the 
description of the $L$-neighbour relation $\sim _L$, is  
coordinate free, see Definition \ref{crucdef} below, and therefore, 
too, is the description of $\Delta f$ and of the notion of harmonic 
function.  We get a characterization of harmonic functions, in terms 
of an average-value property, which is infinitesimal in character and 
does not involve integration, see Theorem \ref{Delta} and Proposition 
\ref{seven}.

\medskip

In Section 3 we prove that diffeomorphisms which preserve the 
$L$-neigh\-bour relation are precisely the conformal ones.  Section 4 
deals with the special case of the complex plane, and Section 5 
explains the ``support of the Laplacian'' in systematic algebraic 
terms.

\section{Preliminaries}
Although the notions we use are introduced in a coordinate 
free way, we have no intention of avoiding use of coordinates 
as a tool of proof. This Section contains in fact mainly 
certain coordinate calculations, which we believe will be 
useful also in other contexts where Riemannian geometry is 
treated in the present synthetic manner.
 
Working in coordinates in $M$ means that we are identifying (an open 
subset of) $M$ with (an open subset of) $R^n$; for simplicity, we talk 
about thse open subsets as if they were all of $M$ and $R^n$, 
respectively; all our considerations are anyway only {\em local}.  The 
Riemannian metric $g$ on $M$ then becomes identified with a Riemannian 
metric on $R^n$, likewise denoted $g$, and it may be written (for 
$x\sim _2 y$) in the form of a matrix product,
$$g(x,y)= (y-x)^T \cdot G(x)\cdot  (y-x),$$
where $x-y\in R^n$ is viewed as a column matrix, and $G(x)$, 
for each $x$, is a symmetric positive definite $n\times n$ matrix.

Using coordinates, we may form affine combinations of (the coordinate 
sets of) points of $M$, at least for sufficiently nearby points, and 
such combinations will in general have only little geometric 
significance, since they depend on the choice of the coordinate 
system.  However, we have the following useful fact:

\begin{prop} Assume $y_1 \sim _1 x$ and $y_2 \sim _1 x$ (so 
$(x+ y_2 - y_1 ) 
\sim _2 x $).  Then
$$g(x, x+ y_2 - y_1 ) = g(y_1 , y_2 );$$
in particular, for $x=0$,
$$ g(0, y_2 - y_1 ) = g(y_1 , y_2 ).$$
\label{simple}\end{prop}

{\bf Proof.} We may assume $x=0$.  Then $g(0, y_2 - y_1 )$ and $g(y_1 
,y_2 )$ are given, respectively, by
$$(y_2 - y_1 )^T \cdot G(0) \cdot (y_2 - y_1 ) = -2 y_1 ^T \cdot G(0) \cdot y_2$$
and
$$(y_2 - y_1 )^T \cdot G(y_1 ) \cdot (y_2 - y_1 ) = -2 y_1 ^T \cdot 
G(y_1 ) \cdot y_2 .$$
Now expand $G(y ) $ as $G(0) + H(y )$ where $H$ depends linearily on 
$y \sim _1 0$.  The difference between our two expressions is then 
$-2 y_1 ^T \cdot H(y_1 ) \cdot y_2$, which depends bilinearily on 
$y_1$ and therefore vanishes.

\medskip

We shall see below (Proposition \ref{good}) that the ``coordinatewise'' affine 
combination considered in Proposition \ref{simple} {\em does } have an 
invariant geometric meaning, provided the coordinate system is {\em 
geodesic}:

We say that 
the metric $g$ on $R^n$ (or equivalently, the coordinate 
system around $x_0 \in M$) is {\em geodesic} at $0\in R^n$ (or 
at $x_0 \in M$, respectively), if the first partial derivatives of 
$G(x)$, as  functions of $x\in R^n$, vanish at $0$; 
equivalently, if $G(x)= G(0)$ for every $x\sim _1 0$. (This is 
in turn equivalent to the vanishing at $x_0$ of the Christoffel symbols 
of the metric, in the given  coordinate system.) It is 
classical that for every point $x_0$, there exists a 
coordinate system which is geodesic at $x_0$. If  
$G(0) $ is the identity matrix, one talks about a  
{\em normal} coordinate system at $x$, and such also exist. 
Cf.\ e.g.\ \cite{CCL} for such notions.

Recall from \cite{Levi-Civita} formula (2) that any Riemannian metric $g: M_{(2)} 
\to R$ admits a unique symmetric extension $\overline{g} : M_{(3)} \to R$;
in coordinates it is given by
\begin{equation} \overline{g} (x,y) = (y-x)^T \cdot ( G(x) + 1/2
(D_{(y-x)}G)(x))\cdot (y-x).\label{corr}\end{equation}

 Recall also from \cite{Levi-Civita} Theorem 3.6 that for $x\sim _2 z$ in a 
 Riemannian manifold, and for $t\in R$, there exists a unique $y_0$ 
 with $y_0 \sim _2 x$ and $y_0 \sim _2 z$ which is a critical point for the 
 function of $y$ given by  \begin{equation}t\overline{g} (x,y) + (1-t) 
 \overline{g} (z,y);
\label{crit}\end{equation}

We call this $y_0$ an (intrinsic) {\em affine combination} of $x$ and $z$.  We 
write it $tx + (1-t)z$; this raises a compatibility problem in case we 
are working in coordinates, since we can then also form the 
``algebraic'' affine combination of two coordinate $n$-tuples.  
However, in geodesic coordinates at $x$, there is no problem, 
according to the following Proposition, which extends Proposition 3.7 
in \cite{Levi-Civita}.  Let us consider a coordinate system with $x$ 
identified with $0$.

\begin{prop} The critical point $y_0$ for the function in 
(\ref{crit}) is the algebraic affine combination $tx + 
(1-t)z$, if either $x\sim _1 z$, or if the coordinate system is 
geodesic at $x$.
\label{alg}\end{prop}

{\bf Proof.} Since $x$ is identified with $0$ in the 
coordinate system, the affine combination in question is just 
$(1-t)z$. To show that it is a critical value for (\ref{crit}) 
means that  
\begin{equation}t\overline{g}(0,(1-t)z +v ) +(1-t) 
\overline{g}(z,(1-t)z +v)\label{hh}\end{equation}
 is independent of $v\sim _1 
0$.  We write $g$ in terms of the symmetric matrices $G$, as above.  
Let us take a Taylor expansion of the function $G(y)$, writing
$$G(y) = G(0) + H(y),$$
where the entries of the matrix $H(y)$ are of degree $\geq 1$ 
in $y$; and if the coordinate system is geodesic at $x=0$, 
$H(y)$ is even of degree $\geq 2$ in $y$. We then calculate.
We get a ``significant'' part from each of the two terms in (\ref{hh}), and then 
some ``error'' terms, each of which will vanish for degree reasons, as 
we shall argue.

The two significant terms are the two terms in

$$t (((1-t)z +v)^T \cdot G(0)\cdot  ((1-t)z +v)) + 
(1-t)( (-tz +v)^T \cdot G(0) \cdot (-tz+v) ).$$

Expanding out by bilinearity and symmetry, the terms involving 
$v$ linearly cancel each other; and the terms involving $v$ quadratically 
vanish because $v\sim _1 0$.  So the significant terms, jointly, do 
not depend on $v\sim _1 0$.

The ``error'' terms are of two kinds: partly, arising from the 
replacement of $G(z)$ by $G(0)$; here, $H(z)$ enters; and partly there are 
correction terms when passing from $g$ to $\overline{g}$ defined on 
pairs of third order neighbours.  The error term of the first kind is 
a multiple of
$$(-tz+v)^T \cdot H(z)\cdot  (-tz+v);$$
we expand this out by bilinearity, and use that $H(z)$ is 
of degree $\geq 1$, and $v\sim _1 0$.  We get four terms each 
of  which vanish for 
degree reasons if {\em either} $z\sim _1 0$ {\em or} if $H(z)$ is of degree 
$\geq 2$.

Finally, the correction  terms for upgrading $g$ to 
$\overline{g}$ don't occur if $z\sim _1 0$, since then $g$ is 
only applied to pairs of second order neighbours. Thus the 
assertion of the Proposition about the case $z\sim _1 x$ is 
already proved. In general, the upgrading involves first 
partial derivatives of $G$, (see (\ref{corr})), so in the 
case the coordinate system is geodesic at $0$, no correction term is 
needed for $\overline{g}(0,(1-t)z)$, but only for $\overline{g}(z, 
(1-t)z+v)$.  Using the formula (\ref{corr}), we see that the correction needed 
is a certain multiple of
$$(-tz + v)^T \cdot ( D_{-tz+v}G)(z) \cdot (-tz+v),$$
hence a linear combination of terms 
$$z\cdot D_z G(z)\cdot z,\; v\cdot D_z G(z)\cdot z,\; z \cdot D_v G(z) \cdot z,$$
and something that contains $v$ in a bilinear way.  All these terms 
vanish for degree reasons: for, since $H$ vanishes in the first 
neighbourhood of $0$, $D_v G(z)$ is of degree $\geq 1$ in $z$, and 
$D_z G(z)$ is even of degree $\geq 2$ in $z$.

\medskip

Essentially the same degree counting as in  this 
proof gives the following result: 
\begin{lemma}Let $y\sim _1 x$ and $z\sim 
_2 x$; then using a geodesic coordinate system at $x=0$, the 
quantity
$\overline{g}(y,z)$ may be calculated as $(z-y)^T \cdot G(0) 
\cdot (z-y)$.
\label{yz}\end{lemma} 

\medskip

Given a Riemannian manifold. If $x\sim _2 z$, the {\em mirror 
image} $z'$ of $z$ in $x$ is by definition the affine 
combination $2x - y$, i.e.\ the $y$ which is critical value for 
$2\overline{g}(x,y) - \overline{g}(z,y)$, 
\cite{Levi-Civita} Theorem 3.6. Also, the parallelogram 
formation $\lambda$ is descibed in \cite{Levi-Civita}. 
Finally, if $t$ is a tangent vector $D\to M$, its geodesic 
prolongation $\overline{t}:D_2 \to M$ is determined by the validity, for 
all $d_1 , d_2 \in D$ of $$\overline{t}(d_1 + d_2 ) = \lambda (t(0), t(d_1 ), 
t(d_2
)).$$
(Recall that $D\subseteq R$ are the elements of square zero, $D_2$ the 
elements of cube zero.) Now Proposition \ref{alg} has the following Corollary:

\begin{prop}
Let $x\sim _2 z$; then the mirror image $z'$ of $y$ w.r.to $x$ 
may be calculated as follows: take  a geodesic 
coordinate system at $x$ with $x=0$. Then $z' = -z$.

Let $y\sim _1 x$, $z\sim _1 x$. Then $\lambda (x,y,z)$ may be 
calculated as follows: take  a geodesic 
coordinate system at $x$ with $x=0$. Then $\lambda (x,y,z) = 
y+z$.

Let $t$ be a tangent vector $D\to M$ at $x\in M$. Then the 
geodesic prolongation $\overline{t} :D_2 \to M$ of $t$ may be 
calculated as follows: take  a geodesic 
coordinate system at $x$ with $x=0$. Let $u$ be the unique 
vector in $R^n$ so that $t(d)= d\cdot u$ for all $d\in D$.  Then for 
$\delta \in D_2$, $\overline{t}(\delta) = \delta \cdot u$
\label{good}\end{prop}

(The vector $u \in R^n$ appearing in the last clause is usually called 
the {\em principal part} of $t$, relative to the coordinate system.)

If $t$ and $s$ are tangent vectors at the same point $x$ of a 
Riemannian manifold $M,g$, we define their inner product 
 $<t,s>$ 
by the validity, for all $d_1 d_2 \in D$, of \begin{equation}d_1 d_2 \;<t,s> = 
-\frac{1}{2} g(t(d_1 ), s(d_2 )).
\label{inner}\end{equation}
In this way, the tangent 
vector space $T_x M$ is made into an inner product space, (and this is 
the contact point with the classical formulation of Riemannian metric).

If $u$ and $v$ are the principal parts of tangent vectors $t$ and $s$ 
at $x\in M$, in some coordinate system at $x=0$ (not necessarily geodesic), 
one has \newline $<t,s> = u^T \cdot G(0) \cdot v$; this follows easily from 
Proposition \ref{simple}.

Combining Lemma \ref{yz} and Proposition \ref{good}, one gets

\begin{lemma} Let $t$ be a tangent vector.  Then for $d\in D$, $\delta 
\in D_2$, we have $$\overline{g} (t(d), \overline{t}(\delta )) 
=(\delta ^2 -
2d\delta )\cdot <t,t>.$$
\label{geo1}\end{lemma}

We are going to define the {\em orthogonal projection} of $z$ 
($z\sim _2 x$) onto a proper tangent $t$ at $x$.  We first define the {\em 
scalar component} of $z$ along $t$; this is unique number $\alpha (z, 
t)$ so that 
\begin{equation}
d\cdot \alpha (z,t) = \frac{1}{2} \frac{(g(x,z)- 
\overline{g}(t(d),z))}{<t,t>}\label{sca}\end{equation} for all $d\in 
D$.  Note that if $z =x$, $\alpha (z,t)=0$, and from this follows that 
for any $z\sim _2 x$, $\alpha (z,t) \sim _2 0$, in other words $\alpha 
(z,t)\in D_2$.  From Lemma \ref{geo1}, applied twice (once with 
$d=0$, once with a general $d\in D$), it is immediate to deduce that if 
$z$ is of the form $\overline{t}(\delta )$ for a $\delta \in D_2$, 
then $\alpha (z, t) = \delta $.

We define the orthogonal projection $proj_{t} (z)$ by $$proj_{t} 
(z)=\overline{t} (\alpha (z,t)).$$ Note that it is a second-order 
neighbour of $x$.  It follows  from the above that if $z$ 
is of the form $\overline{t}(\delta )$, then $proj _t (z)=z$.  

 \section{Laplacian neighbours}

Here is the crucial definition:

\begin{defn}Let $z\sim _2 x$. We say that $z$ is a {\em  Laplacian 
neighbour}  of $x$ (written $z\sim _L x$) if for 
every proper tangent $t$ at $x$, we have \begin{equation}g(x,z)= n 
\cdot g(x, proj_t (z)),\label{cruceqn}\end{equation} where $n$ is the 
dimension of the manifold.
\label{crucdef}\end{defn}

Maybe one of the names ``isotropic, harmonic, or conformal, neighbour" 
would be more appropriate.

\medskip

Clearly $z\sim _1 x$ implies $z\sim _L x$; for if $z$ is a 
first-order neighbour of $x$, then so is its 
orthogonal projection, and hence both the $g$-quantities to be 
compared in (\ref{cruceqn}) are zero.  If the dimension $n$ is 1, $\sim _L$ 
is the same as $\sim _2$ but in general, the set ${\cal M}_{L} (x)$ of 
$L$-neighbours of $x $ is much smaller than the set ${\cal M}_2 (x)$ of 
second-order neighbours; in fact, the ring of functions on ${\cal M}_{L} (x)$ 
is a finite dimensional vector space which is just one dimension 
bigger than the ring of functions on ${\cal M}_{1} (x)$, as we shall see in 
the proof of Proposition \ref{uni-c} below.

  We conjecture that the relation $\sim _L$ is symmetric, but we haven't 
  been able to do the necessary calculations, except in the 
  case of $R^n$, where the symmetry is easy to prove, using 
  Proposition \ref{five} below.

Note the following curious phenomenon in dimension $n\geq 2$:
 if $z\sim _L x$, then $z$ does not connect to $x$ by 
{\em any} geodesic $D_2\to M$ (given by a proper tangent vector $t$), 
except perhaps in the trivial case when  $g(x,z)=0$.  In 
other words, the $L$-neighbours of $x$ are genuinely {\em isotropic}, in the 
sense that they are in {\em no} preferred direction $t$ (hence the 
alternative name ``isotropic neighbour'' suggested).  Nevertheless, there 
are sufficiently many $L$-neighbours of $x$ to define the Laplacian 
differential operator $\Delta$, see Theorem \ref{Delta} below.

  Let us assume the manifold 
in question has dimension $n$.  Then we have

\begin{prop}In any geodesic normal coordinate system at $x=0$, $z=(z_1 , 
\ldots ,z_n )$ is $\sim _L 0$ if and only if
 $$z_i ^2 = z_j ^2 \mbox{ for all $i,j$, and } z_i z_j =0 \mbox{ for }i\neq j$$
 (and $z_i ^3 =0$ for all $i$; this latter condition follows from the 
 other two if $n\geq 2$).
 \label{five}\end{prop}
{\bf Proof.} First, if $t$ and $s$ are tangent vectors at $x=0$ with principal 
parts $u$ and $v$, respectively (meaning $t(d) = du, 
s(d)=dv$), then 
$<t,s>= u\bullet v$, where $\bullet$ denotes the usual dot product of 
vectors in $R^n$.  Also, if $t$ is a tangent at $x=0$ with principal 
part $u$, then $\alpha (z,t ) = (z\bullet u)/(u\bullet u)$; for, 
calculating the enumerator in (\ref{sca}) gives (using Proposition 
\ref{simple}) $$z\bullet z - \overline{g}(du,z)= z\bullet z - 
(z-du)\bullet (z-du) = 2d\; z\bullet u.$$ From the third clause in 
Proposition \ref{good}, we then get the familiar looking 
\begin{equation}proj _t (z) = \frac{z\bullet u}{u\bullet u} u.
\label{familiar}\end{equation}

In particular, if $t$ is the (proper) tangent vector with principal part $e_i 
\in R^n$ ($= (0, \ldots, 1, \ldots 0)$ (with 1 in the $i$'th position, 
$0$'s elsewhere), then $proj _t (z_1 , \ldots ,z_n ) = z_i e_i$.  In 
particular $g(0, proj _t (z)) = z_i ^2$.  If, on the other hand, $t$ 
is the tangent vector with principal part $e_{i,j}$ (the vector 
with $1$'s in the $i$'th and in the $j$'th position, $i\neq j$, $0$'s 
elsewhere), then $proj _t (z)$ has $(z_i +z_j)/2$ in the $i$'th and in 
the $j$'th position, and $0$'s elsewhere.  In particular, $$g(0, proj 
_t (z))=\frac{1}{2}(z_i ^2 + z_j ^2) + z_i z_j .$$ If $z$ therefore is 
an $L$-neighbour of $0$, we conclude that $z_i ^2 = z_j ^2$ for all 
$i,j$, and that $z_i z_j =0$ if $i\neq j$.

Conversely, assume that in some geodesic normal coordinate system 
at $x=0$, the coordinates $(z_1 , \ldots ,z_n )$ satisfy the equations 
$z_i ^2 = z_j ^2$, $z_i z_j =0$ for $i\neq j$, and let $t$ be a proper 
tangent vector at $x$ with principal part $u =(u_1 , \ldots ,u_n )$.  
Then
 $$proj _t (z) = \frac{z\bullet u}{u\bullet u} u,$$
 and therefore 
 $$g(x,proj _t (z)) = (\frac{z\bullet u}{u\bullet u} u)\bullet 
( \frac{z\bullet u}{u\bullet u} u),$$
which we calculate by arithmetic to be
$$\frac{(\sum _i u_i z_i )(\sum _j u_j z_j )}{u \bullet u}= \frac{\sum 
_{ij } u_i u_j z_i z_j}{u\bullet u},$$
but since $z_i z_j =0$ for $i\neq j$, only the ``diagonal'' terms 
survive, and we are left with
$$\frac{\sum _i u_i u_i z_i z_i}{u\bullet u}.$$
But $z_i z_i = z_1 z_1$ for all $i$, so this factor can go outside the 
sum sign in the enumerator, and we 
get $z_1 ^2 (\sum_i u_i u_i)/u\bullet u = z_1 ^2$, which is $1/n$ 
times $\sum z_i ^2$ since all the $z_i ^2$ are equal.  This proves the 
Proposition.

\medskip

From  Propositions \ref{good} and \ref{five}, one 
immediately deduces that if $z\sim _L x$, then also $z'\sim _L x$ for 
any affine combination $z' = t x + (1-t )z$ ($t \in R$).

\begin{prop}If two functions $f_1$ and $f_2$: ${\cal M}_L (x)\to R$ agree on 
${\cal M}_1 (x)$ there is a unique number $c\in R$ such that for all $z\sim 
_L x$
$$f_1 (z)-f_2 (z) =c\cdot g(x,z).$$
\label{uni-c}\end{prop}

{\bf Proof.} 
 Using a geodesic normal 
coordinate system at $x=0$, it is a matter of analyzing the ring of functions
${\cal M}_L (0) \to R$ for the case where $M= R^n$ with standard 
inner-product metric.  The Proposition gives that ${\cal M}_L (0)$ may be 
described as $D_L (n) \subseteq R^n$, defined by 
\begin{equation}D_L (n):= \{ (d_1 , \ldots , d_n )\in R^n \mid d_i ^2 = d_j ^2
\mbox{, and } d_i d_j =0 \mbox{ for } i\neq j \} 
,\label{DLN}\end{equation} (for $n\geq 2$; for $n=1$, $D_L (1) = D_2 = 
\{\delta \in R \mid \delta ^3 = 0\}$).  This is (for $n\geq 2$) the 
object represented by the Weil algebra ${\cal O}(D_L (n)):=k[Z_1 , 
\ldots , Z_n ]/I$, where $I$ is the ideal generated by the $Z_i ^2 - 
Z_j ^2$, and by $Z_i Z_j$ for $i\neq j$.  It is immediate to calculate 
that, as a vector space, this ring is $(n+2)$-dimensional, with
linear generators $$1, Z_1 , \ldots ,Z_n , Z_1 ^2 + \ldots + Z_n ^2 .$$  
By the general (Kock-Lawvere) axiom scheme for SDG \cite{SDG}, \cite{lav}, 
this means that any function $f:D_L (n) \to R$ is of the form
$$f(z_1 , \ldots , z_n ) = a + \sum _i b_i z_i + c (\sum _i z_i ^2 ),$$
for unique $a, b_1 , \ldots ,b_n , c \in R$, or equivalently
$$f(z_1 , \ldots , z_n ) = a + \sum _i b_i z_i + c g(0,z).$$
Since the restriction of $f$ to $D(n)$ is given by the data $a, b_1 , 
\ldots ,b_n$, the unique existence of $c$ follows.

\medskip

The following Theorem deals with an arbitrary Riemannian manifold $M, 
g$ of dimension $n$, and gives a coordinate free characterization of the 
Laplacian operator $\Delta$.

\begin{thm}For any $f:{\cal M}_L (x)\to R$, there is a unique number $L$ 
with the property that for any $z\sim _L x$
$$f(z) + f(z') -2f(x) = L\cdot g(x,z),$$
where $z'$ denotes the mirror image of $z$ in $x$. We write 
$\Delta f (x) := n L$.

\label{Delta}\end{thm}

Put differently,
$$f(z) + f(z') -2f(x) = \frac{\Delta f (x)}{n} g(x,z).$$
If the function $f$ is {\em harmonic} at $x$, meaning that $\Delta f 
(x) =0$, it follows that it has a strong {\em average value property}: 
the value at $x$ equals the average value of $f$ over any pair of 
points $z$ and $z'$ ($L$-neighbours of $x$) which are symmetrically 
located around $x$.

\medskip

{\bf Proof.} Again, we pick a normal geodesic coordinate system at 
$x=0$, so identify ${\cal M}_L (x)$ with $D_L (n)$; then $z'$ gets identified 
with $-z$, by Proposition \ref{good}.  The left hand side of the 
expression in the Theorem then has restriction 0 to $D(n)$, being 
(with notation as above) $(a + \sum b_i z_i )+(a + \sum b_i (-z_i )) - 
2a.$ Hence the unique existence of $L$ follows from 
Proposition \ref{uni-c}.

\medskip

The following Proposition serves to as partial justification of the
 use of the name ``Laplacian'' for the $\Delta $ considered in the 
Theorem. We consider the standard Riemannian metric on $R^n$, 
$g(x,z)= ||z-x||^2$, for $z\sim _2 x$.

\begin{prop}Let $f:R^n \to R$ and let $x\in R^n$. Then
$$\Delta f (x) = \sum _i \frac{\partial ^2 f}{\partial x_i ^2} (x).$$
\label{six}\end{prop}

{\bf Proof.}  For simplicity, let $x=0$, so that $z' = -z$. 
We Taylor expand $f(z)$ and $f(-z)$ 
from $0$, and consider $f(z) + f(-z) -2f(0)$; then terms of degree 
$\leq 1$ cancel, and we get
$$\sum _{ij}\frac{\partial ^2 f}{\partial x_i \partial x_j}z_i z_j + 
\mbox{ higher terms};$$ 
 now the calculation proceeds much like the 
one in the proof of Proposition \ref{five} above: if $z\sim _L 0$, 
only the diagonal terms in the sum survive, all the $z_i ^2$ are equal 
to $z_1 ^2$, which we move outside the parenthesis, and get $z_1 ^2 $ 
times the classical Laplacian $\sum \partial ^2 f/ \partial x_i ^2 
(0)$.  But $z_1 ^2 = 1/n \cdot g(0,z)$.

\medskip

Similarly, one proves by Taylor expansion

\begin{prop} If $z\sim _L x$ in $R^n$, then for any $f:R^n \to 
R$,
$$f(z)= f(x) + df_x (z-x) + \frac{1}{2n} \Delta f (x) 
||z-x||^2 .$$ 
\label{Taylor} \end{prop}

Recall that any function $M\to R$ looks affine on any 1-monad ${\cal 
M}_1 (x)$; functions that look affine on the larger $L$-monads ${\cal 
M}_L (x)$ are precisely the harmonic ones:

\begin{prop}Assume $f:M\to R$ is harmonic at $x$.  Then for an $z\sim 
_L x$, $f$ preserves affine combinations of $x$ and $z$. 
Conversely, if for given $x$, $f$ preserves affine 
combinations of $x$ and $z$ for every $z\sim _L x$, then $f$ 
is harmonic at $x$.  (``Harmonic" at $x$ here in the sense: $\Delta f (x)=0$.)
\label{seven}\end{prop}
(Recall that the affine combination $tx + (1-t)z$ is defined as the 
critical point $y_0$ in (\ref{crit}).)

\medskip

{\bf Proof.} Assume $z\sim _L x$, and pick a geodesic normal 
coordinate system at $x=0$. Without loss of generality, we may 
assume $f(0)=0$. Then to say that $f$ preserves 
affine combinations of $x$ and $z$ is to say that for all 
$s\in R$, $f(sz)=sf(z)$.
For $z \sim _L 0$, we have by Proposition \ref{Taylor} that $f(z) = \sum a_i z_i 
+ c\sum z_i ^2$ for unique $a_i$ and $c$ ($c= \Delta f (0) /2n$); $f(sz)$ and $sf(z)$ 
have the same terms of first order in $z$; their second order terms 
are respectively $c s^2 \sum z_i ^2$ and $s c \sum z_i ^2$, and if 
these two expressions are to be equal for all $s$ and all $z$, we must 
have $c=0$.  This means that $f$ is harmonic at $x$.  Conversely, if 
$f$ preserves affine combinations of the kind mentioned, it preserves 
the affine combination $2x -z$, or equivalently, the left hand side of 
the expression in Theorem \ref{Delta} is 0, hence it follows that $L$ 
and hence $\Delta f (x) $ is 0.

\medskip

{\bf Remark.} With some hesitation, I propose to call a map $f:M \to N$ between 
Riemannian manifolds {\em harmonic} if it preserves affine 
combinations of $L$-neighbours in $M$ and if it preserves the property 
of being $L$-neighbours.  I have not been able to compare the proposed 
definition, with a certain classical concept of harmonic map between 
Riemannian manifolds.  But at least: When the codomain is $R$, the 
definition is the basic classical one of harmonic function, by 
Proposition \ref{seven}.  For, preservation of the $\sim _L$ relation 
is automatic when the codomain is $R$, since in $R$, $\sim _L$ is the same 
as $\sim _2$.

\section{Conformal maps}

We consider a diffeomorphism $f:M\to N$ between Riemannian 
manifolds $(M, g), (N, h)$.  
To say that $f$ is an {\em isometry} at $x\in M$ is to say that for all $z\sim 
_2 x$,  $g(x,z)=h(f(x),f(z))$.  To say that $f$ is {\em conformal} at 
$x\in M$ with constant $k = k(x) > 0$ is to say  that for all 
$z\sim _2 x$, $h(f(x),f(z)) = k(x) g(x,z)$ (so if $k(x)=1$, $f$ is an 
isometry at $x$).  The terminology agrees with classical usage, as we 
shall see below.  We first prove

\begin{prop} Assume $f:M \to N$ is conformal at $x\in M$ with 
\newline $h(f(x),f(z)) = k g(x,z)$.  Then for all $y_1 \sim _1 x$, $y_2 \sim 
_1 x$,
$$h(f(y_1 ), f(y_2 )) = k g(y_1 ,y_2 ),$$
and conversely.
\label{conf}\end{prop}

{\bf Proof.} We choose coordinates, and assume $x=0$ and $f(x) =0$; the 
metrics in $M$ and $N$ are then given by functions $g$ and $h$, 
respectively, and they are in turn given by symmetric matrices $G(y)$, 
and $H(z)$ for all $y\in M$ and $z\in N$.  We now calculate $k g(y_1 
,y_2 )$.  We have, by Proposition \ref{simple} that $$k g(y_1 , y_2 )= 
k g (0, y_2 - y_1 ) =
h(0, f(y_2 -y_1 )).$$
Now there is a bilinear $B(-,-)$ such that for all pairs of 
$1$-neighbours $y_1 , y_2$ of $0$, we have $f(y_2 - y_1) = f(y_2 ) - 
f(y_1 ) + B(y_1 ,y_2 )$.  So the calculation continues $$= h(0 , f(y 
_1 ) - f( y_2 ) + B(y_1 ,y_2 ))$$ $$ = (f(y _1 ) - f( y_2 ) + B(y_1 ,y_2 
))^T \cdot H(0) \cdot (f(y _1 ) - f( y_2 ) + B(y_1 ,y_2 )).$$ Since 
$f$ depends in a linear way of $y_1 \sim _1 0$ and $y_2 \sim _1 0$, this whole 
expression multiplies out by linearity, and for degree reasons all terms 
involving $B$, as well as some others, vanish, and we are left with 
$-2 f(y_1 )^T \cdot H(0) \cdot f(y_2 )$.  On the other hand, $h(f(y _1 
), h(f(y_2 )) = h(0, f(y_2 ) - f(y _1 ))$, by Proposition 
\ref{simple}, and writing this in terms of $H(0)$ gives the same 
expression.

The converse is proved in the same way for $z\sim _2 0$ of the form 
$y_2 - y_1 $ with $y_1 \sim _1 0$ and $y_2 \sim _ 1 0$, but this 
suffices to get the result for all $z\sim _2 0$, by general principles 
of SDG (``$ R$, and hence any manifold, perceives the addition map 
$D(n)\times D(n) \to D_2 (n)$ to be epic''.)

\medskip

Call a linear map $F:U\to V$ between inner product vector spaces {\em 
conformal} with constant $k > 0$ if for all $u_1 ,u_2 \in U$
$$<F(u _1 ),F(u_2 )> = k <u_1 , u_2 >.$$
It follows immediately from Proposition \ref{conf}, and from the 
construction of inner product in the vector space of tangents at $x$, 
and at $f(x)$, that if $f$ is conformal at $x$ with constant $k$, then 
$ df_x : T_x M \to T_{f(x)} N$ is a conformal linear map with the 
same constant $k$.  The converse also holds; for if $df _x$ is 
conformal with constant $k$, we deduce that for all pairs of tangents 
$t$ and $s$ at $x$
$$g(f(t(d_1 )), f(s(d_2 ))) = k \cdot g(t(d_1 ), s(d_2 )),$$
and hence 
\begin{equation}g(f(y_1 )),f(y_2 )) = k \cdot g(y_1 , y_2 
)\label{conv}\end{equation}
 for all $y_i$'s of the form $t(d)$ for a 
tangent vector $t$ and a $d\in D$.  Again by general principles, any 
manifold $N$ ``perceives all 1-neighbours of $x$ to be of this form''.  
From Proposition \ref{conf} we therefore deduce that $f$ is conformal 
at $x$ with constant $k$.

\begin{thm}A diffeomorphism $f$ is conformal at $x\in M$ if and only if 
$f$ maps ${\cal M}_L (x)$ into ${\cal M}_L (f(x))$.
\label{three}\end{thm}

{\bf Proof.} Assume $f$ maps ${\cal M}_L (x) $ into ${\cal M}_L (f(x))$. We may 
pick normal coordinates at $x$ as well as at $f(x)$.  The 
neighbourhoods ${\cal M}_L (x)$ and ${\cal M}_L (f(x))$ then both get 
identified with $D_L (n)$, and $x=0$, $f(x)=0$.  The restriction of 
$f$ to $D_2 (n)$, $f:D_2 (n) \to R^n$, takes $0$ to $0$ and is therefore 
of the form $f(y) = A\cdot y + B(y)$, where $A$ is an $n\times n$ 
matrix, and $B(y)$ is a map $R^n \to R^n$ which is homogeneous of 
degree 2 in $y\in R^n$, i.e.\ an $n$-tuple of quadratic forms $B_i$.

Assume now that $f$ maps $D_L (n)$ into itself.  For $z\in D_L (n)$, 
the $i$'th coordinate of $f(z)$ is
$$f_i (z) = \sum _k a_{ik} z_k + B_i (z).$$
Squaring this, only the terms in $( \sum _k a_{ik} z_k)( \sum _l 
a_{il} z_l)$ survive for degree reasons (using that $z\in D_2 (n)$).  
But using further that $z_k z_l =0$ for $k\neq l$, only the ``diagonal'' 
terms survive, and we get 
\begin{equation}f_i (z ) ^2 =\sum _k a_{ik} ^2 z_k ^2 = z_1 ^2 \sum _k 
a_{ik} ^2 .\label{sq}\end{equation} Similarly for $i\neq j$
 \begin{equation}f_i (z) f_j (z)=z_1 ^2 ( \sum _k a_{ik} 
 a_{jk}).\label{ij}\end{equation}
If now $f(z) \in D_L (n)$ for all $z\in D_L (n)$, we get that the 
expression in (\ref{sq}) is independent of $i$, and from the 
uniqueness assertion in Proposition \ref{uni-c} we therefore conclude
$$ \sum _k a_{ik}^2 = \sum _k a_{jk}^2 \mbox{ for all }i,j ;$$
and similarly we conclude from (\ref{ij}) that
$$\sum _k a_{ik} a_{jk} =0 \mbox{ for } i\neq j.$$
These two equations express that all the rows of the matrix $A$ have the 
same square norm $k$, and that they are mutually orthogonal.  This 
implies that the linear map $df_x$ represented by the matrix is 
conformal, and hence $f$ is conformal at $x$.  

The proof that conformality of $f$ at $x$ implies that $f$ maps ${\cal M}_L 
(x)$ into ${\cal M}_L (f(x))$ goes essentially through the same calculation, 
and is omitted.
 
 \section{A famous pseudogroup in dimension 2}
 
 The content of the present section is partly classical, namely the 
 equivalence of the various ways of describing the notion of 
 holomorphic map from (a region in) the complex plane $C=R^2$ to 
 itself.  Synthetic concepts enter essentially in two of the conditions 
 in the Theorem below, namely 1) and 7).
 
 \medskip

An {\em almost complex} structure on a general manifold $M$ 
consists in giving, for each $x\in M$, a map $I_x :{\cal M}_1 (x)\to 
{\cal M}_1 (x)$ with $I_x (x)= x$ and $I_x (I_x (z)) = z'$ for any 
$z\sim _1 x$; Here, $z'$ denotes the mirror image of $z$ in 
$x$, i.e.\ the affine combination $2x -z$; recall \cite{DFIC} 
that affine combinations of 1-neighbours make ``absolutely'' sense, 
i.e.\ do not depend on, say, a Riemannian structure.  It is clear what 
it means for a map $f$ to {\em preserve} such structure at the point 
$x$: $f(I_x (z)) = I_{f(x)} (f(z))$.

The manifold $R^2$ carries a canonical almost-complex 
structure, given by
$$I_{(x_1 ,x_2 )} (z_1 , z_2 ) = (x_1 - (z_2 - x_2 ), x_2 + 
(z_1 - x_1 )).$$
Identifying $R^2$ with the complex plane $C$, this is just 
$$I_x (z) = x + i(z-x).$$
Utilizing the multiplication of the complex plane $C$, we may 
consider the set $D_C$ of elements of square zero in $C$ (recalling 
the fundamental role which the set $D $ of elements of square 
zero in $R$ plays in SDG). We have, by trivial calculation,

\begin{prop} Under the identification of $C$ with $R^2$,
$$D_C = D_L (2).$$
\label{eleven}\end{prop}
Having $D_C$, we may mimick the basics of SDG and declare a 
function $f: C\to C$ to be {\em complex differentiable} at 
$x\in C$ if there is a number $f'(x)\in C$ so that
$$f(z) = f(x) + f'(x)\cdot (z-x) \mbox{ for all }z \mbox{ with } 
z-x \in D_C .$$
(The uniqueness of such $f'(x)$, justifying the notation, 
follows from the general axiom scheme of SDG, applied to 
$D_1 (2)$, the 1-jet classifier in $R^2$.  Note $D_1 (2) \subseteq 
D_L (2)$.)  --- The notion of course makes sense for functions $f$ which are 
just defined {\em locally} around $x$.

\begin{thm}Let $f:R^2 \to R^2$ be a (local) orientation 
preserving diffeomorphism.  Let $x\in R^2 =C$.  Then the following conditions 
are equivalent:

1) $f$ maps ${\cal M}_L (x)$ into ${\cal M}_L (f(x))$

2) $f$ is conformal at $x$

3) $f$ satisfies Cauchy-Riemann equations at $x$

4) $f$ preserves almost complex structure at $x$.

\noindent Also the following conditions are equivalent, and they imply 1)-4):

5) $f$ is complex-differentiable at $x$ 

6) $f$ maps ${\cal M}_L (x)$ into ${\cal M}_L (f(x))$, and $f$ preserves affine 
combinations of $x$ and $z$ for any $z\sim _L x$.  

\noindent Finally, if 1)-4) hold for all $x$, 5) and 6) hold for all $x$.
\label{famous}\end{thm} 
(Note that 6) says that $f$ is {\em harmonic} 
at $x$, in the sense of Remark at the end of Section 2.)

\medskip

{\bf Proof.} The equivalence of 1) and 2) is already in 
Theorem \ref{three}, and this in turn is, as we have seen, equivalent to 
conformality of the linear $df _x$.  But conformal orientation 
preserving $2\times 2$ matrices are of the form
\begin{equation}\left[ \begin{array}{rr}a& -b\\b & a \end{array} 
\right] .\label{CR}\end{equation}
 Since the entries of the matrix for $df _x$ are 
$\partial f_i / \partial x_j$, this form (\ref{CR}) of the matrix 
therefore expresses that the Cauchy-Riemann equations hold at 
$x$, i.e. is equivalent to 3). On the other hand, a simple 
calculation with $2\times 2$ matrices give that a matrix 
commutes with the matrix $I = \left[ \begin{array}{rr} 
0&-1\\1&0 \end{array} \right]$  for the almost complex 
structure iff it has the above ``Cauchy-Riemann" form (\ref{CR}).

Now assume 5). If $f$ is complex differentiable at $x$, we prove that condition 1) 
holds at $x$ as follows.  Let $z\sim _L x$.  Then $(z-x)^2 =0$ 
Proposition \ref{eleven}, and by complex differentiability
\begin{equation}f(z)-f(x)= f'(x)(z-x),\label{co}\end{equation}
 so since the right hand side has 
square zero, then so does the left hand side, but again by Proposition 
\ref{eleven}, this means that $f(z)\sim _L f(x)$,  proving 1), 
and hence also the first part of 6).  But also,  if $f$ is 
complex-differentiable at $x$, $f$ preserves affine 
combinations of the form $tx + (1-t)z$ for $z\in {\cal M}_L (x)$; this 
follows from (\ref{co}), since $z-x \in D_C (x)$ by Proposition 
\ref{eleven}, so also the second part of 6) is proved.  Conversely, if 
6) holds, $f$ is conformal at $x$ by Theorem \ref{three}, so $df_x$ is 
if the form (\ref{CR}).  Then $f'(x) = a+ib$ will serve as the complex 
derivative; for, since $f$ preserves affine combinations of $x$ and 
$z$ for $z\sim _L x$, we have the first equality sign in
$$f(z)= f(x)+ df_x (z-x) = f(x) + f'(x)\cdot (z-x)$$ 
for such $z$, i.e.  for $z-x \in D_L (2)= D_C $.

Finally, assume 1)-4) hold for all $x$.  Then we may differentiate the Cauchy-Riemann 
equations for $f = (f_1 , f_2 )$   by 
$\partial /\partial x_1$ and $\partial / \partial x_2$
and compare, arriving in the standard way to $\Delta f_1 \equiv 0$ and 
$\Delta f_2 \equiv 0$.  From the ``Taylor expansion'' in Proposition 
\ref{Taylor}, applied to $f_1$, we conclude that $f_1 (z) = f_1 (x) + 
(df_1 )_x (z-x)$ for $z\sim _L x$, and similarly for $f_2$, so $f(z)= 
f(x)+ df_x (z-x)$ for such $z$, i.e.\ for $z-x \in D_C$.  Since $df_x$ 
is given by a conformal matrix (\ref{CR}), by 2), this proves that 
$a+ib$ will serve as the complex derivative of $f$ at $x$.

\section{Support of the Laplacian}
We arrived at $D_L (n)$ from the geometric side, namely as the 
$\sim _L$-neighbours of $0$ in the Riemannian manifold 
$M=R^n$; the differential operator $\Delta$ was then seen to provide 
the top term in the Taylor expansion of functions defined on 
$D_L (n)$.

Here, we briefly indicate how to arrive at $D_L (n)$ from the 
algebraic side, starting with $\Delta = \sum \partial ^2  /\partial 
x_i ^2$. More precisely, we consider $\Delta$ as a {\em 
distribution} at $0\in R^n$. So $\Delta$ is the linear map
$$k[X_1 , ... ,X_n ] \to R$$
given by
\begin{equation}f\mapsto \sum _i \frac{ \partial ^2 f} {\partial 
x_i ^2}(0).
\label{star}\end{equation}
The algebraic concept that will give $D_L (n)$ out of this 
data is the notion of {\em coalgebra}, and {\em subcoalgebra}, as in 
\cite{swe}.  If we let $A$ denote the algebra $k[X_1 , ...  ,X_n ]$, 
then the distribution $\Delta$ of (\ref{star}) factors
$$A\to B \to k,$$
where $A\to B$ is an {\em algebra} map, and $B$ is finite 
dimensional (take e.g.\ $B=A/J$ where $J$ is the ideal 
generated by monomials of degree $\geq 3$). The set $A^{o}$ of 
linear maps $A\to k$ having such a factorization property 
constitute a coalgebra, \cite{swe} Proposition 6.0.2. Every 
element in a coalgebra generates a finite dimensional 
subcoalgebra, by \cite{swe} Theorem 2.2.1. In particular, 
$\Delta \in A^{o}$ generates a finite dimensional coalgebra 
$[\Delta ]$ of $A^{o}$, and this coalgebra {\em ``is"} $D_L 
(n)$. More specifically, the dual algebra of $[\Delta ]$ is 
the coordinate ring ${\cal O} (D_L (n))$ of $D_L (n)$, i.e.\ the Weil 
algebra ${\cal O} (D_L (n)):= k[Z_1 , ...  ,Z_n ]/I$ considered in the 
proof of Proposition \ref{uni-c}, as we shall now argue.

The following ``Leibniz rule" for $\Delta$ is well known,
$$\Delta (f\cdot g) = \Delta f \cdot g + 2\sum _i \frac{\partial 
f}{\partial x_i} \cdot \frac{\partial g}{\partial x_i} + 
f\cdot \Delta g .$$
This means that in the coalgebra $A^{o}$, we have the 
following formula for $\psi (\Delta )$ ($\psi =$ the 
comultiplication of the coalgebra; $\delta $ the Dirac distribution 
``evaluate at 0''): \begin{equation}\psi (\Delta ) = \Delta \otimes 
\delta  +2\sum _i \frac{\partial }{\partial x_i} \otimes \frac{\partial 
}{\partial x_i} + \delta \otimes \Delta ,
\label{twostar}\end{equation}
where now $\Delta$, $\partial d /\partial x_i$, $\delta$ are viewed as 
distributions at $0$, like in (\ref{star}), meaning that one evaluates 
in $0$ after application, $$f\mapsto (\Delta f)(0)\mbox{, }f\mapsto 
\frac{\partial
f}{\partial x_i}(0)\mbox{, } f\mapsto f(0).$$
From (\ref{twostar}), (and from $\psi (\partial /\partial x_i 
)=\partial /\partial x_i \otimes \delta + \delta \otimes 
\partial /\partial x_i$, which expresses the Leibniz rule for 
$\partial /\partial x_i$)
we see that the subcoalgebra $[\Delta ]$ 
generated by $\Delta$ is generated as a vector space by the 
elements
$$\delta \mbox{, }\frac{\partial}{\partial x_1} , ...
,\frac{\partial}{\partial x_n}\mbox{, } \Delta ,$$
and since these are clearly linearly independent, we see that $[\Delta ] 
\subseteq A^{o}$ is $(n+2)$-dimensional.  The dual algebra of $[\Delta 
]$ is a quotient algebra $A/I$ of $A$, where $I$ is the ideal of those 
$f\in A$ which are annihilated by the elements of $[\Delta ]$.  This 
ideal $I$ contains $x_i ^2 - x_j ^2$, and $x_i x_j$ for $i\neq j$.  
Since the quotient of $A$ by the ideal generated by $x_i ^2 - x_j ^2$, 
and $x_i x_j$ for $i\neq j$ is already $(n+2)$-dimensional, as 
calculated in the proof of Proposition \ref{uni-c}, it follows that 
the quotient algebra there is actually the dual of $[\Delta ]$.

\medskip

The idea that a coalgebra like $[\Delta ]$ {\em is} itself an 
infinitesimal geometric object goes back to Gavin Wraith in 
the early seventies, \cite{wraith}.  The specific way of generating Weil 
algebras from differential operators was considered by Emsalem 
\cite{Emsalem} (without coalgebras).


\begin{thebibliography}{99} 

\bibitem{CCL}S.S.\ Chern, W.H.\ Chen and K.S.\ Lam, Lectures on 
differential geometry, World Scientific 1999.

\bibitem{Emsalem} J.\ Emsalem, G\'{e}om\'{e}trie des points \'{e}pais, 
Bull.\ Soc.\ math.\ France 106 (1978), 399-416.

\bibitem{FM} A.\ Kock, Formal manifolds and synthetic theory 
of jet bundles, Cahiers de Topologie et G\'{e}om\'{e}trie 
Diff\'{e}rentielle 21 (1980), 227-246.

\bibitem{SDG} A.\ Kock, Synthetic Differential Geometry, London Math.\  
Soc.\ Lecture Notes Series No.\ 51, Cambridge Univ.\ Press 1981.

 \bibitem{Levi-Civita} A.\ Kock, 
Geometric construction of the Levi-Civita Parallelism, Theory and 
Appl.\ of Categories 4 (1998), No.\ 9.

\bibitem{VF} A.\ Kock, Volume form as volume of infinitesimal 
simplices, arXiv:math.CT/0006008, June 2000.

\bibitem{DFIC} A.\ Kock, Differential forms as infinitesimal 
cochains, Journ.\ Pure Appl.\ Alg.\, to appear.

\bibitem{lav} R.\ Lavendhomme, Basic concepts of Synthetic 
Differential Geometry, Kluwer Texts in the Math.\ Sciences 13, Kluwer 
1996.

\bibitem{swe} M.E.\ Sweedler, Hopf Algebras, Benjamin 1969.

\bibitem{wraith} G.C.\ Wraith, Talk at Oberwolfach, July 27, 1972.

\end{thebibliography}
\end{document}